\documentclass{amsart}
\usepackage{setspace,amsthm,amsmath,amscd,latexsym}
\usepackage[psamsfonts]{amssymb}
\usepackage[cmtip,arrow,matrix,curve]{xy}
\usepackage{pb-diagram,pb-xy}

\newtheorem{thm}{Theorem}[section]
\newtheorem{lem}[thm]{Lemma}
\newtheorem{prop}[thm]{Proposition}
\newtheorem{cor}[thm]{Corollary}

\theoremstyle{definition}
\newtheorem{defi}[thm]{Definition}

\theoremstyle{remark}
\newtheorem{rmk}[thm]{Remark}
\newtheorem{ex}[thm]{Example}

\newcommand\Om{\Omega}
\newcommand \n{\nabla}
\newcommand\del{\partial}
\newcommand\vp{\varphi}

\newcommand\tom{\widetilde \Omega}

\newcommand{\sdr}[2]{\overset{#1}{\underset{#2}{\rightleftharpoons}}}
\newcommand{\Z}{\mathbb{Z}}

\begin{document}

\title{A chain coalgebra model for the James map}
\author{Kathryn Hess}
\author{Paul-Eug\`ene Parent}
\author{Jonathan Scott}

\address{\'Ecole Polytechnique F\'ed\'erale de Lausanne \\
  Institut de G\'eom\'etrie, Alg\`ebre et Topologie \\
  CH-1015 Lausanne \\
  Switzerland}
\address{Department of Mathematics and Statistics \\
    University of Ottawa \\
    585 King Edward Avenue \\
    Ottawa, ON \\
    K1N 6N5 Canada}
\address{Department of Mathematics and Statistics \\
    University of Ottawa \\
    585 King Edward Avenue \\
    Ottawa, ON \\
    K1N 6N5 Canada}

    \email{kathryn.hess@epfl.ch}
    \email{pparent@uottawa.ca}
    \email{jscott@uottawa.ca}

\date{\today}

\maketitle

\section{Introduction}

Let $L$ be a $1$-reduced simplicial set.  Let $G(-)$ and $C(-)$ be the Kan loop group
and normalized chain functors respectively. The explicit,
natural twisting cochain $t_{L} : CL
\rightarrow C(GL)$ of Szczarba \cite{S} determines a natural morphism of
chain algebras $\theta_{L} : \Omega CL \rightarrow C(GL)$
that induces an isomorphism in homology, since $CL$
is $1$-connected.  Here $\Omega CL$ is the cobar construction
\cite{A}.

The coassociative, counital diagonal $\psi$
on $\Omega CL$ of Baues \cite{B} makes $\Omega CL$ a Hopf algebra. The
Alexander-Whitney diagonal on $CL$ is comultiplicative up to strong
homotopy or DCSH \cite{GM}, and hence is the linear part in a morphism of chain
algebras, $\tilde{\Omega}\Delta \colon \Omega(CK) \rightarrow
\Omega(CL \otimes CL)$. By \cite{HPST},
$\psi=q\circ\tilde{\Omega}\Delta$, where $q : \Omega (CL \otimes CL) \rightarrow \Omega
CL \otimes \Omega CL$ is the Milgram equivalence \cite{Mi}.\footnote{In hindsight, we
realize that we had resurrected an idea of Drachman \cite{D1,D2}.}
Furthermore, the Szczarba equivalence is a
DCSH morphism.

In the present paper, we consider the special case $L = EK$, where
$K$ is a pointed simplicial set and $E(-)$ is the simplicial
suspension.  Note that we allow the case when $K$ is not reduced,
and hence $L$ need not be $1$-reduced.  Since the Alexander-Whitney
diagonal on $CEK$ is trivial, there is a natural isomorphism of
chain algebras,
\begin{equation}\label{eq:iso}
    \Omega CEK \cong T(\tilde{C}K),
\end{equation}
where $\tilde{C}K = CK / C(*)$. Our main result states:

\medskip

\noindent{\bf Theorem} (Theorem \ref{thm:b-s-canon})
The isomorphism (\ref{eq:iso}) preserves diagonals.

\medskip

To prove the theorem, we construct a chain model of the James map
(see below), $\alpha_{K} : CK \rightarrow \Omega CEK$. We compute
the diagonal $\psi_{EK}$ on $\Omega CEK$, even when $EK$ is not
$1$-reduced, and show that $\alpha_{K}$ is comultiplicative.  An
immediate consequence is that $\psi_{EK}$ is coassociative.

Let $JX$ be the James construction \cite{J} on the pointed topological space $X$.
The unit $\eta_{X} : X \rightarrow JX$ determines an
isomorphism $T\tilde{H}_{*}(X) \xrightarrow{\cong} H_{*}(JX)$ of associative algebras
when
the coefficients are such that $\tilde{H}_{*}(X)$ is torsion-free.
This result has the following implications at the chain level.  Let
$CS(-)$ be the normalized singular chains functor over $\Z$ (the
composite of normalized chains $C$ and total singular complex $S$).
Since $JX$ is a topological monoid, $CS(JX)$ is a chain Hopf
algebra.  The chain coalgebra morphism $CS(\eta_{X}) : CS(X) \rightarrow
CS(JX)$, extends to a chain Hopf
algebra morphism, $T(CS(X)) \rightarrow CS(JX)$ inducing an isomorphism in
homology.

If one replaces $X$ by a pointed simplicial set $K$, our main
theorem identifies $T(CK) \cong \Omega CEK$ \emph{as a Hopf
algebra}. We calculate the Szczarba equivalence $\Omega CEK
\rightarrow C(GEK)$ when $K$ is reduced, and show that it factors
through the natural equivalence $T(CK) \rightarrow C(G^{+}EK)$ where
$G^{+}EK$ is a simplicial model of $J|K|$.
We are interested in $G^+EK$ since in the applications,
it is important to be able to treat $\Omega
CEK$ when $K$ is not necessarily reduced.

We remark that Theorem~\ref{thm:b-s-canon} is used in an essential way in
\cite{H}. In future work, we will apply the ideas in the present
article to construct a sequence of operads that control iterated
cobar constructions.

In the Appendix, we show that the Milgram map $q:\Omega(A\otimes C)\rightarrow\Omega
A\otimes\Omega C$ fits into a strong
deformation retract of chain complexes \cite{GM}, and therefore is
a chain homotopy equivalence whenever $A$ and $C$ are coaugmented. This result
further reinforces the geometric validity of our calculation when
$K$ is not reduced.

\section{Notation and background}

\subsection{Simplicial structures}

For any $m\leq n\in\mathbb N$, let
$[m,n]=\{j\in\mathbb N\vert m\leq j\leq n\}$. Let $\boldsymbol{\Delta}$ denote the
category with objects $Ob(\boldsymbol{\Delta})=\{[0,n]\vert n\geq0\}$
and morphisms
$$
\boldsymbol{\Delta}([0,m],[0,n])=\{[0,m]\stackrel{f}{\rightarrow}[0,n]\,\vert\, f\hbox{ is
an order-preserving set map}\}.
$$
In \cite{Mac} p.177, it is shown that the classical
{\it coface} and {\it codegeneracy} maps,
namely,

\begin{minipage}{4.5cm}
\begin{eqnarray*}
d_i^n&:&[0,n-1]\rightarrow[0,n]\\
x&\mapsto&
\left\{
\begin{array}{cc}
x&\hbox{ if $x<i$}\\
x+1&\hbox{ if $x\geq i$}
\end{array}\right.
\end{eqnarray*}
\end{minipage}
\hspace{1cm}and\hspace{1cm}
\begin{minipage}{4.5cm}
\begin{eqnarray*}
\sigma_i^n&:&[0,n+1]\rightarrow[0,n]\\
x&\mapsto&
\left\{
\begin{array}{cc}
x&\hbox{ if $x\leq i$}\\
x-1&\hbox{ if $x>i$}
\end{array}\right.
\end{eqnarray*}
\end{minipage}

\bigskip

\noindent for $0\leq i\leq n$,
generate $\boldsymbol{\Delta}([0,m],[0,n])$.
Eilenberg and MacLane in \cite{EM1}
associate to a morphism
$f\in\boldsymbol{\Delta}([0,m],[0,n])$ its {\it derived}
map $f'\in\boldsymbol{\Delta}([0,m+1],[0,n+1])$ defined by
$$
f'(0)=0\quad\hbox{and}\quad f'(i+1)=f(i)+1,\quad
i=0,...,m.
$$
Clearly, $(d_i^n)'=d_{i+1}^{n+1}$ and
$(\sigma_i^n)'=\sigma_{i+1}^{n+1}$.

A {\it simplicial set} is a contravariant functor $K:\boldsymbol{\Delta}\rightarrow Set$.
Let $K_n:=K([0,n])$,
$s^n_i:=K(\sigma_i^n)$, and $\partial^n_i:=K(d^n_i)$. The
maps $s^n_i$ and $\partial^n_i$ are called respectively {\it degeneracy}
and {\it face} maps and an element $x\in K_n$ is called an
$n$\emph{-simplex}. The dimension superscript will be omitted when the context
is clear. We will use extensively the simplicial
identities:
\begin{eqnarray*}
\partial_i\partial_j&=&\partial_{j-1}\partial_i\quad\hbox{if}\quad
i<j\hbox{,}\\
s_is_j&=&s_{j+1}s_i\quad\hbox{if}\quad i\leq j\hbox{,}\\
\partial_is_j&=&s_{j-1}\partial_i\quad\hbox{if}\quad i<j\hbox{,}\\
\partial_js_j&=&\hbox{identity}\quad=\quad\partial_{j+1}s_j\hbox{, and}\\
\partial_is_j&=&s_j\partial_{i-1}\quad\hbox{if}\quad i>j+1.
\end{eqnarray*}

The derived maps associated to the face and degeneracy maps are
$$
(\partial_i^n)'=K((d_i^n)')=\partial_{i+1}^{n+1}\quad\hbox{and}\quad
(s_j^n)'=K((\sigma_j^n)')=s_{j+1}^{n+1}.
$$
Naturality implies that the derived map associated to a
composition of iterated faces and degeneracies is the composition
of the iterated derived face and degeneracy maps, i.e.,
$(s_{i_1}...s_{i_k}\partial_{j_1}...\partial_{j_l})'=
s_{i_1+1}...s_{i_k+1}\partial_{j_1+1}...\partial_{j_l+1}$.
Given two simplicial sets $K$ and $L$, we extend the notion of
derived map on $K\times L$ componentwise, i.e.,
\begin{eqnarray*}
\big((s_{i_1}...s_{i_k}\partial_{j_1}...\partial_{j_l})\times
(s_{\alpha_1}...s_{\alpha_n}\partial_{\beta_1}...\partial_{\beta_m})\big)'&=&\\
(s_{i_1}...s_{i_k}\partial_{j_1}...\partial_{j_l})'&\times&
(s_{\alpha_1}...s_{\alpha_n}\partial_{\beta_1}...\partial_{\beta_m})'.
\end{eqnarray*}

Let $\Delta^n$ denote the {\it standard geometric $n$-simplex}.
The {\it singular complex} on a topological space $X$ is the simplicial set $S(X)$ where
$S_n(X)=Top(\Delta^n,X)$. Its left adjoint, the {\it geometric realization}
functor,
is denoted $|\cdot|$.

More generally a simplicial object in a category $\mathcal C$
is a contravariant
functors from $\boldsymbol{\Delta}$ to $\mathcal C$. In particular
we will be concerned with simplicial monoids and simplicial
(abelian) groups.

Let $\mathcal F_{ab}:Set\rightarrow Ab$ denote the free abelian
group functor. Given a simplicial set $K$, there is an associated simplicial abelian group,
namely $K_{ab}:=\mathcal F_{ab}\circ K$. We extend linearly the notion of
derived map.
For all $n>0$, let $DK_n=\cup_{i=0}^{n-1}s_i(K_{n-1})$, the set
of degenerate $n$-simplices of $K$. The {\it normalized chain complex} on $K$, denoted
$C(K)$, is given by
$$
C_n(K)=\mathcal F_{ab}(K_n)/\mathcal F_{ab}(DK_n),
$$
where its differential is induced by $\partial=\sum_{i=0}^n(-1)^i\partial_i$.
The simplicial identities imply that
$\partial^2=0$ and $\partial(\mathcal F_{ab}(DK_n))\subset\mathcal
F_{ab}(DK_{n-1})$. As noted in \cite{EM1} the notion of derived
map does not descend to normalized chains. Hence all computations
will be carried out with unnormalized chains, and we pass
afterwards to the quotient.

Let $*$ be the unique simplicial set generated by a single
non-degenerate $0$-simplex. For the remainder of this paper we work within
the category of {\it pointed} simplicial sets.
Objects are simplicial maps $\kappa:*\rightarrow K$, which boils down to choosing
a $0$-simplex as basepoint in $K$. Morphisms are commuting
triangles, i.e., basepoint-preserving simplicial
maps.
The basepoint of $K$ is denoted $k_0:=\kappa(*)$, and
$k_n:=(s_0)^nk_0$. An $n${\it-reduced} simplicial set $K$ is a simplicial set such
that $K_i=\{k_i\}$ for $0\leq i\leq n$. A $0$-reduced simplicial set will
simply be called a \emph{reduced} simplicial set. The unique simplicial map
$c:K\rightarrow*$ induces a chain map
$C(K)\stackrel{\epsilon}{\longrightarrow}
C(*)\cong\mathbb Z$. The {\it reduced chain complex} on $K$ will be
denoted by
$
\widetilde{C}(K)=\hbox{ker}(\epsilon)
$.
Thus $\widetilde{C}_{>0}(K)=C_{>0}(K)$, while
$\kappa:*\rightarrow K$ induces a
natural splitting of $\epsilon$, from which we obtain
a natural basis for $\widetilde{C}_0$, namely the set
$
\big\{x-k_0\,|\,x\in K_0\backslash\{k_0\}\big\}
$.
Thus $C_0(K)=\mathbb Z\{k_0\}\oplus\widetilde{C}_0(K)$.

We now recall three classical constructions on a pointed simplicial set
$K$ that are central to the discussion at hand.

\noindent{\bf(a) Simplicial suspension} $EK$ (\cite{M} p.124): Let
$E_0(K)=b_0$, while $E_n(K)$ is the
set of pairs $(i,x)$, where $i \geq 1$ is an integer and $x \in
K_{n-i}$, under the identification $(i,k_n)=s_0^{n+i}b_0=b_{n+i}$.
The face and degeneracy operators are generated by
\begin{enumerate}
    \item $\partial_0(1,x)=b_{n}$, for all $x \in K_{n}$,
    \item $\partial_1(1,x)=b_0$, for all $x\in K_0$,
    \item $\partial_{i+1}(1,x)=(1,\partial_ix)$, for all $x \in
    K_{n}$, $n>0$,
    \item $s_0(i,x) = (i+1,x)$, and
    \item $s_{i+1}(1,x)=(1,s_ix)$,
\end{enumerate}
with all other face and degeneracy maps defined by the requirement
that $EK$ be a simplicial set.

The reduced suspension of the
geometric realization of $K$ is canonically homeomorphic to the
geometric realization of $EK$ (\cite{M} p.125). Note that $C_{n+1}(EK)$
is generated by elements $(1,x)$, for nondegenerate $x \in K_{n}$.

\noindent{\bf(b) Simplicial loop group} $G(K)$ (\cite{M} p.118): Let $K$
be a reduced simplicial set. Define $G_n(K)$ to be the free group
generated by the elements of $K_{n+1}$ under the identification
$s_0x=e_n$ for $x\in K_n$, where $e_n$ is the identity element of
$G_n(K)$. The face and degeneracy operators are given by
\begin{enumerate}
\item $\tau(\partial_0x)\partial_0\tau(x)=\tau(\partial_1x)$,
\item $\partial_i\tau(x)=\tau(\partial_{i+1}x)$, if $i>0$, and
\item $s_i\tau(x)=\tau(s_{i+1}x)$, if $i\geq0$,
\end{enumerate}
where $\tau(x)$ denotes the class of $x\in K_{n+1}$ in $G_n(K)$.
These maps, $\partial_i$ and $s_i$, extend uniquely to
homomorphisms $G_n(K)\rightarrow G_{n-1}(K)$ and
$G_n(K)\rightarrow G_{n+1}(K)$ respectively.

In the Appendix of
\cite{Sm}, Smith shows that the (based) loops on the
geometric realization of $K$ is weakly equivalent to the
geometric realization of $GK$.

\noindent{\bf(c) Simplicial James construction} $G^+E(K)$:
Notice that since $G(K)$ is a simplicial group, from condition (b.1)
we can solve
for $\partial_0\tau(x)$. If we ask that $G_n(K)$ be the free
monoid instead of the free group, then condition (b.1) is not
enough to determine $\partial_0\tau(x)$. But if $K=EL$, then
condition (b.1) characterizes $\partial_0\tau(x)$. Indeed,
for $x\in L_{n>0}$, we have either
$$\tau(\partial_0(1,x))=\tau(b_n)=e_{n-1}\quad\hbox{or}\quad
\partial_0\tau(i+1,x)=\partial_0\tau(s_0(i,x))=e_{n+i-1}.
$$
Hence we have a well defined functor, $G^+E(-)$, from (arbitrary)
pointed simplicial sets to simplicial monoids.
An easy calculation shows that for all $K$, the natural map
\begin{eqnarray*}
\eta _{K}:K&\rightarrow& G^+E(K)\\
x&\mapsto&\tau {(1,x)}
\end{eqnarray*}
is simplicial.  In fact, using the universal properties of the James construction, $JX$,
on a topological space $X$ (\cite{J}) and the adjunction between
$S(-)$ and $|\cdot|$, one can show that $\eta_K$ is a model of the topological James map
$X\rightarrow JX$. In particular, $|G^+E(K)|\cong J|K|$. Moreover,
when $K$ is reduced, the inclusion $G^+E(K)
\hookrightarrow GE(K)$ is a homotopy equivalence. In
\cite{M} p.126, the last result is said to be valid for countable reduced simplicial sets.
Using Proposition 2.4 on p.9 of
\cite{GJ} one can extend it to arbitrary reduced simplicial sets.

\begin{rmk}
The $G^+E(-)$ construction is isomorphic to Milnor's $F^+$
construction as remarked in \cite{Sm}.
\end{rmk}

\subsection{Differential structures}

We recall now a number of basic definitions and constructions related to graded modules and graded (co)algebras over a
principal ideal domain $R$.
A graded $R$-module $V=\oplus_{i\in\mathbb
Z}V_i$ is {\it connected} if $V_{<0}=0$ and $V_0\cong R$. It is
{\it simply connected} if, in addition, $V_1=0$. We write $V_+$
for $V_{>0}$.
Let $V$ be a non-negatively graded, free $R$-module. The free associative
algebra generated by $V$ is denoted $TV$, i.e.,
$$
TV\cong R\oplus V\oplus(V\otimes V)\oplus(V\otimes V\otimes
V)\oplus...
$$
where the product $\mu:TV\otimes TV\rightarrow  TV$
is given by word concatenation. We denote
the submodule of words of length $n$ by $T^nV=V^{\otimes n}$.

The {\it suspension} endofunctor $s$ on the category of graded
modules is defined on objects $V=\oplus_{i\in\mathbb
Z}V_i$ by $(sV)_i\cong V_{i-1}$. Given a homogeneous element $v\in
V$, we write $sv$ for the corresponding element of $sV$. The
suspension $s$ admits an obvious inverse, which we denote
$s^{-1}$. Observe that $C(EK) \cong sC(K)$ as chain complexes.
A map of chain complexes inducing an isomorphism in homology will
be called a {\it quasi-isomorphism}.

 Let $f,g:(A,d)\rightarrow(B,d)$ be two maps of
chain algebras.  An \emph {$(f,g)$-derivation} is a linear map $\vp
:A\rightarrow B$ of degree $+1$ such that  $\vp \mu
=\mu (\vp \otimes g+f \otimes \vp)$, where $\mu$ denotes the
multiplication on $A$ and $B$.
A \emph{derivation homotopy} from $f$ to $g$ is an $(f,g)$-derivation
$\vp $ that satisfies $d\varphi + \varphi d = f - g$.

Let $(C,d,\Delta)$ be a coaugmented chain coalgebra. Let $\overline{C}=
\ker\varepsilon$ where $\varepsilon:C\rightarrow R$ is the counit. The {\it reduced
coproduct} is defined by
$$\overline{\Delta}(c):=\Delta(c)-(c\otimes1+1\otimes c),$$
for $c\in \overline{C}$.
\begin{defi}[\cite{A}] The cobar construction on $C$,
denoted $\Omega(C)$, is the chain algebra
$(Ts^{-1}(\overline{C}),d_{\Omega})$, where
$d_{\Omega}=-s^{-1}ds+\mu(s^{-1}\otimes s^{-1})\overline{\Delta}s$
on generators.
\end{defi}
A word in $T^n(s^{-1}(\overline{C}))$ will be denoted by
$[x_1\vert...\vert x_n]:=s^{-1}x_1\otimes...\otimes s^{-1}x_n$,
while the unit will be denoted by $[\,]$. There is a natural chain
algebra morphism
\begin{eqnarray}\label{Milgram}
q:\Omega(C\otimes C')\rightarrow\Omega(C)\otimes\Omega(C')
\end{eqnarray}
specified by $q([x\otimes1])=[x]\otimes[\,]$, $q([1\otimes
y])=[\,]\otimes[y]$, and $q([x\otimes y])=0$ for all $x\in \overline{C}$ and
$y\in \overline{C}'$. Milgram shows in \cite[Theorem 7.4]{Mi}
that if $C$ and $C'$ are 1-connected, then $q$ is a natural
quasi-isomorphism of chain algebras. In Appendix A we extend this result to
arbitrary coaugmented
chain coalgebras.

Let $(C,d,\Delta)$ be a chain coalgebra, and let
$(A,d,\mu)$ be a chain algebra.  A \emph{twisting cochain} from $C$
to $A$ is a degree $-1$ map $t:C\rightarrow  A$ of graded modules such
that
$$
dt+td=\mu (t\otimes t)\Delta.
$$
If $C$ is connected, then any twisting cochain $t:C\rightarrow  A$
induces a chain algebra map $\theta :\Omega (C)\rightarrow A$ by setting
$\theta([c])=t(c)$. It is equally clear that any chain algebra map
$\theta :\Omega (C)\rightarrow A$ gives rise to a twisting cochain via the
composition
$$C_{+}\stackrel {s^{-1}}{\longrightarrow} s^{-1}  C_{+}\hookrightarrow Ts^{-1}  C_{+}\stackrel{\theta}{\longrightarrow} A.$$

In section 4 of this paper, we work in the category {\bf DCSH} \cite{GM}.
Its objects are augmented, connected coassociative chain
coalgebras. A DCSH-morphism from $C$ to $C'$ is a map of chain algebras
$\Omega(C)\rightarrow  \Omega(C')$.
In a slight abuse of terminology, we say that a chain map between chain
coalgebras $f:C\rightarrow  C'$ is a DCSH-\emph{map}
if there is a morphism  in $\mathbf{DCSH}(C,C')$ of which $f$ is the linear part.
In other words, for $c\in \overline{C}$, there is a map of chain algebras
$g:\Omega(C)\rightarrow \Omega(C')$ such that
$$g([c])=[f(c)]+\text{higher-order terms}.$$
In a further abuse of notation, we sometimes
write $\widetilde \Omega f:\Omega C\rightarrow \Omega C'$
to indicate one choice of chain algebra map of which $f$ is the linear part.

\subsection {Homological perturbation theory}

We now recall those elements of homological perturbation
theory that we need for this article.

\begin{defi} Suppose that $\n :(X, \del ) \rightarrow (Y,d)$ and
$f:(Y,d)\rightarrow (X,\del )$ are morphisms of chain
complexes.  If $f\n = 1_{X}¥$  and there exists a
chain homotopy
$\vp : (Y,d)\rightarrow (Y,d)$ such that
\begin{enumerate}
\item $d\vp +\vp d =\n f -1_{Y}¥$,
\item $\vp \n =0$,
\item $f\vp =0$, and
\item $\vp \sp 2=0$,
\end{enumerate}
then $(X,d) \sdr{\n}f (Y,d)\circlearrowleft\vp$ is a \emph{strong
deformation retract (SDR) of chain complexes.}
\end{defi}

The following notion was introduced by Gugenheim and Munkholm.

\begin{defi} An SDR $(X,d) \sdr{\n}f (Y,d)\circlearrowleft\vp$ is
called \emph{Eilenberg-Zilber (E-Z) data} if $(Y,d, \Delta _{Y})$ and
$(X,d, \Delta _{X})$ are chain coalgebras and $\n$ is  a morphism of
coalgebras.
\end{defi}

Observe that in this case
$$(d\otimes 1_{X}¥+1_{X}¥\otimes d)\bigl((f\otimes f)\Delta _{Y}¥ \vp\bigr)+\bigl((f\otimes
f)\Delta_{Y}¥ \vp\bigr)
d=\Delta_{X}¥ f-(f\otimes f)\Delta_{Y}¥,$$
i.e., $f$ is a map
of coalgebras up to chain homotopy. In fact, as the following theorem of
Gugenheim and Munkholm shows, $f$ is usually a DCSH map.

\begin{thm}\label{thm:g-m} \cite[Theorem 4.1]{GM} Let $(X,d) \sdr{\n}f
(Y,d)\circlearrowleft\vp$ be E-Z data such that $X$ is simply
connected and $Y$ is connected. Let $F_0=0$ and $F_{1}$ be the composite
$$
Y\stackrel{f}{\longrightarrow}X\rightarrow
X_+\stackrel{s^{-1}}{\longrightarrow}s^{-1}X_+,
$$
and construct inductively $F_{k}:Y\rightarrow T^k(s^{-1}X_+)$ by the formula
$$
F_{k}=-\sum _{i+j=k}(F_{i}\otimes F_{j})\Delta _{Y}¥ \vp.
$$
Then $\displaystyle F=\prod_{k\geq1}F_k=\bigoplus_{k\geq1}F_k$ is a twisting cochain.
Similarly, let $\Phi_0$ be the natural augmentation on $Y$ and $\Phi _{1}$
be the composite
$$
Y\stackrel{\varphi}{\longrightarrow}\underbrace{Y\rightarrow
Y_+\stackrel{s^{-1}}{\longrightarrow}s^{-1}Y_+}_{\rho_Y},
$$
and construct inductively $\Phi _{k}:Y\rightarrow T^k(s^{-1}Y_+)$ by the formula
$$
\Phi _{k}=\left(\Phi_{k-1}\otimes \rho_Y + \sum _{i+j=k}\Bigl
(\Omega(\n)F_j\otimes\Phi _{i}\Bigr )\right)\Delta _{Y}¥ \vp.
$$
Then $\displaystyle\Phi=\prod_{k\geq0}\Phi_k=\bigoplus_{k\geq0}\Phi_k$ is a twisting (homotopy) cochain.
Moreover,
$$
\Om(X,d) \sdr{\Om\n}{\tom f}
\Om(Y,d)\circlearrowleft\tom \vp
$$
is an SDR, where $\Omega\nabla$ is the algebra morphism determined
by the coalgebra morphism $\nabla$, $\tom f$ is the algebra
morphism determined by the twisting cochain $F$, and $\tom\vp$ is
the derivation homotopy determined by the twisting cochain $\Phi$.
\end{thm}

\begin{rmk}
Since $X$ is simply connected, the key, as Gugenheim and Munkholm noted,
is that on elements of degree $n$, $F_k=0$
when $k\geq n$. Thus $F$ is a well defined map into
$\bigoplus_{k\geq1}T^k(s^{-1}X_+)$ and not merely into
$\prod_{k\geq1}T^k(s^{-1}X_+)$. That is the only place
where the hypothesis $X$ be simply connected is used. Thus that
condition can be removed, and hence Theorem \ref{thm:g-m} still applies
if for some other reason (geometric,
algebraic,...) $F$ turns out to be locally nilpotent.
\end{rmk}

\subsection{Relating simplicial and differential structures}

Let $K$ and $L$  be simplicial sets. The natural \emph{Alexander-Whitney map} is the morphism of chain complexes $f_{K,L}:C(K\times
L)\rightarrow  C(K)\otimes C(L)$ given on nondegenerate  $x\times y\in (K\times
L)_n$ by
$$
f_{K,L}(x\times y)=\sum_{i=0}^n\tilde{\partial}^{n-i}x\otimes\partial_0^iy,
$$
where
$\tilde{\partial}^{n-i}:=\partial_{i+1}...\partial_n$. The
simplicial diagonal $K\stackrel{\Lambda}{\rightarrow }K\times K$ together with
the Alexander-Whitney map endow the normalized chains on $K$ with
a natural coproduct $\Delta _{K}=f_{K,K}\circ C(\Lambda)$ (\cite{EM1,EM2,M}).  If $x\in K_n$ is nondegenerate, then
$$
\Delta_K(x)=\sum_{i=0}^n\tilde{\partial}^{n-i}x\otimes\partial_0^ix.
$$

It is well known that the Alexander-Whitney map is a chain
equivalence \cite{EM2}. A natural chain homotopy inverse is the {\it Eilenberg-Zilber
map} $\nabla_{K,L}:C(K)\otimes C(L)\rightarrow  C(K\times
L)$ defined by
$$
\nabla_{K,L}(x\otimes
y)=\sum_{(\mu,\nu)}(-1)^{\epsilon(\mu)}s_{\nu_q}...s_{\nu_1}x\,\times\,s_{\mu_p}...s_{\mu_1}y
$$
for nondegenerate $x\in K_p$ and $y\in L_q$, where the sum is taken over
all $(p,q)$-shuffles $(\mu,\nu)$ and
$\epsilon(\mu)=\sum_{i=1}^p[\mu_i-(i-1)]$ is the signature of the
corresponding permutation.
Recall that a
$(p,q)$-shuffle is a permutation $\pi$ of $\{0,...,p+q-1\}$ such
that $\pi(i)<\pi(j)$ if $0\leq i<j\leq p-1$ or $p\leq i<j\leq
p+q-1$. We use Eilenberg and MacLane's convention where we let $\mu_i=\pi(i-1)$,
$1\leq i\leq p$, and $\nu_j=\pi(j+p-1)$, $1\leq j\leq q$. Clearly
$\pi$ is determined by $\mu$ and $\nu$, and we let
$\pi=(\mu,\nu)$.
Notice that, in contrast to the Alexander-Whitney map,
$\nabla_{K,L}$ is a coalgebra map.

There is a natural chain homotopy, $\vp_{K,L}
$ such that
\begin{eqnarray}\label{EMEZ}
C(K)\otimes C(L)\sdr {\n_{K,L} }{f_{K,L}}C(K\times L)\circlearrowleft\vp_{K,L}
\end{eqnarray}
constitutes E-Z data.  If $K$ and $L$ are 1-reduced, then Theorem~\ref{thm:g-m} implies that $f_{K,L}$ is naturally a DCSH map.
In \cite{EM2} Eilenberg and MacLane give a recursive formula for $\vp _{K,L}$
that we reproduce here.
Set $\varphi(C_0(K\times L))=0$. For $q>0$ we let
$$
\varphi(x\times y)\,=\,-\varphi'(x\times y)\,+\,(\nabla_{K,L}\circ
f_{K,L})'s_0(x\times y),
$$
for $x\times y\in (K\times L)_q$ nondegenerate.

\subsection{The cobar diagonal}
We now describe the coproduct structure on $\Omega(C(K))$ for a $1$-reduced simplicial set
$K$.
It was defined in \cite {HPST} and shown there to be identical to the Baues coproduct.
Recall that the Alexander-Whitney map $f_{K,L}$ is naturally a DCSH map. We are thus led to define
\begin{eqnarray*}\label{psik}
\psi _{K}\,=&
q\circ\tilde{\Omega}f_{K,K}\circ\Omega C(\Lambda)&
\end{eqnarray*}
where $q$ is the Milgram equivalence (see (\ref{Milgram}) p.\pageref{Milgram}).
The cobar diagonal, $\psi _{K}$, is strictly coassociative and
cocommutative up to derivation homotopy.
Furthermore, the Szczarba quasi-isomorphism of chain algebras
$\theta _{K}:\Om C(K)\rightarrow C(GK)$ is a DCSH map,
in a way compatible with the algebra structures.

\section{A chain model for the simplicial James map.}

Let $K$ be an arbitrary pointed simplical set.  Recall that $EK$ is reduced,
so the chain coalgebra $C(EK)$ is connected and $\Omega C(EK)$ is defined.
Define a map
\begin{eqnarray}\label{alpha}
\alpha:C(K)\rightarrow\Omega C(EK)
\end{eqnarray}
by $\alpha(k_0)=[\,]$,
$\alpha(y-k_0)=[(1,y)]$ if $y\in K_0\backslash\{k_0\}$, and
$\alpha(x)=[(1,x)]$ if $x\in K_{\geq1}$.

We proceed to show that $\alpha$ is a chain map.
As in the geometric case, the coproduct on $C(EK)$ is as simple as possible.

\begin{lem}\label{lem:primitive-diagonal}
The diagonal on $C(EK)$ is primitive, i.e.,
$\overline{\Delta}_{EK}=0$.
\end{lem}
\begin{proof}
The face map $\partial_{0}$ applied to any element $(1,x) \in
E_{n+1}(K)$  is the basepoint $b_{n}$, and $b_{n}$ is degenerate
unless $n=0$.  Therefore
\begin{eqnarray*}
    \Delta(1,x) & = & \sum_{j=0}^{n+1} \partial_{j+1} \cdots
        \partial_{n+1}(1,x) \otimes \partial_{0}^{j}(1,x) \\
                & = & b_{0} \otimes (1,x) + (1,x) \otimes b_{0}
\end{eqnarray*}
as desired.
\end{proof}

\begin{cor}
The differential in $\Omega C(EK)$ is linear, and
\[
    d_{\Omega}[(1,x)] =
    \sum_{j=0}^{n}(-1)^{j}[(1,\partial_{j}x)],
\]
for $x\in K_n$.
\end{cor}
\begin{proof}
The quadratic part of the differential vanishes by
Lemma~\ref{lem:primitive-diagonal}.  A straightforward computation gives the result.
\end{proof}

\begin{thm}\label{thm:geom-model-adj}
The map $\alpha$ is a chain map.
\end{thm}
\begin{proof}
We only check for $x\in K_1$ since it is obvious in other degrees.
Recall that $(1,k_0)=b_1$ is degenerate in $EK$. We have
\begin{eqnarray*}
\alpha(\partial x)&=&\alpha(\partial_0x-\partial_1x)\\
&=&\alpha(\partial_0x-k_0-(\partial_1x-k_0))\\
&=&\alpha(\partial_0x-k_0)-\alpha(\partial_1x-k_0)\\
&=&\left\{\begin{array}{cc}
[(1,\partial_0x)]-[(1,\partial_1x)]&\hbox{ if }\partial_0x\neq
k_0\hbox{ and }\partial_1x\neq k_0\\
{[(1,\partial_0x)]}&{\hbox{ if }\partial_0x\neq
k_0\hbox{ and }\partial_1x=k_0}\\
{-[(1,\partial_1x)]}&\hbox{ if }\partial_0x=k_0\hbox{ and }\partial_1x\neq k_0\\
0&\hbox{ if }\partial_0x=k_0\hbox{ and }\partial_1x=k_0
\end{array}\right.\\
&=&d_{\Omega}\alpha(x).
\end{eqnarray*}
\end{proof}
The chain map $\alpha$ is in fact a chain model of the simplicial
James map
$K \rightarrow G^+EK$.  More precisely, the diagram below commutes
\begin{eqnarray}\label{triangle}
\xymatrix{
&{\Om C(EK)}\ar[d]^{\gamma}\\
C(K)\ar[ur]^{\alpha}\ar[r]_-{C(\eta_K)}&C(G^+EK),
}
\end{eqnarray}
where the chain algebra map $\gamma$ is induced by the twisting
cochain $(1,x)\mapsto\tau(1,x)$. Hence diagram (\ref{triangle})
induces a commuting chain algebra triangle
\begin{eqnarray}\label{Algtriangle}
\xymatrix{
&{\Om C(EK)}\ar[d]^{\gamma}\\
T\widetilde{C}(K)\ar[ur]^{\widehat{\alpha}}\ar[r]_-{\widehat{C}(\eta_K)}&C(G^+EK),
}
\end{eqnarray}
where $\widehat{\alpha}$ and $\widehat{C}(\eta_K)$ are the chain
algebra morphisms induced by $\alpha$ and $C(\eta_K)$
respectively. Since $C(\eta_K)$ is a coalgebra morphism, $\widehat{C}(\eta_K)$ is a
Hopf algebra morphism. A homological argument shows that
$\widehat{C}(\eta_K)$ is a quasi-isomorphism while
$\widehat{\alpha}$ is obviously an isomorphism. Hence $\gamma$
is a quasi-isomorphism.

\section{Bott-Samelson, Szczarba and the cobar diagonal}

In this section, we extend the definition of the cobar
diagonal to the suspension of an arbitrary pointed simplicial set $K$. Indeed, we
show that we have a factorisation
$$
\xymatrix{
C(EK)\ar@{.>}[rd]\ar[r]^-{C(\Lambda)}&C(EK\times
EK)\ar[r]^-{\prod_{k\geq1}F_k}&\prod_{n\geq0}T^ns^{-1}{(C(EK)\otimes
C(EK))_+}\\
&\Omega(C(EK)\otimes C(EK))\ar@{^{(}->}[ru]\ar[r]^-q&\Omega(C(EK))\otimes\Omega(C(EK))
}
$$
where $\{F_k\}_{k\geq1}$ is the family of higher homotopies
associated to the Alexander-Whitney map $f_{EK,EK}$ as given by
Theorem \ref{thm:g-m}. Moreover, we show that the composite of the dotted
arrow together with $q$, denoted $\xi_{K}$, is a twisting cochain.
We define the (extended) \emph{cobar diagonal} to be
the induced chain algebra map
\begin{eqnarray*}
\psi_{K}:\Omega(C(EK))&{\longrightarrow}&{\Omega(C(EK))\otimes\Omega(C(EK))}\\
{[(1,x)]}&{\mapsto}&\xi_{K}\bigl((1,x)\bigr).
\end{eqnarray*}
Note that when $K$ is reduced, then
$\psi_{K}=q\circ\tilde{\Omega}f_{K,K}\circ\Omega C(\Lambda)$ as on
page \pageref{psik}.

Finally, we prove Theorem \ref{thm:b-s-canon}, a chain-level
Bott-Samelson theorem with respect to the (extended) cobar
diagonal.  We then show that the
Szczarba equivalence is a strict morphism of chain Hopf algebras
for suspensions. In the process we establish a few helpful
combinatorial identities involving face maps in simplicial
suspensions.

\begin{rmk}
 From
Lemma~\ref{lem:primitive-diagonal} we deduce that $\Delta_{EK}$ is
a coalgebra map. Thus one can directly obtain a diagonal
on $\Omega C(EK)$ as $q\circ\Omega(\Delta_{EK})$.
If $\psi_{K}=q\circ\Omega(\Delta_{EK})$ it would
imply, however, that the higher homotopies were trivial. As we show below,
it is not true in general:  the homotopy $F_2$ is usually nonzero.
Endowed with the ``wrong" coproduct,
$\Omega (C(EK))$ is indeed a chain Hopf algebra,
but it is not  weakly equivalent to $C(G^+EK)$.
\end{rmk}

\subsection{The (extended) cobar diagonal}

To show that $F_2$ is nonzero in general, we need first to study
the Eilenberg-MacLane homotopy $\varphi : C(EK \times EL)
\rightarrow  C(EK \times EL)$ (\cite{EM2}). Since we will
eventually be interested only in the image of $C(\Lambda)$, i.e.,
simplices of the form $(1,x)\times(1,x)$ for nondegenerate
$(n-1)$-simplices $x$, we will assume that $n\geq1$ and
concentrate on simplices of the form $(1,x)\times(1,y)\in(EK\times
EL)_n$ for nondegenerate $x\in K_{n-1}$ and $y\in L_{n-1}$.
Recall that since $EK$ and $EL$ are reduced, $C(EK \times EL)$ is
a connected coalgebra and so the reduced diagonal is defined.

\begin{prop}\label{prop:form-diag}
On simplices of the form $(1,x)\times(1,y)$, we have
$$\bar{\Delta}_{EK\times EL}\circ\varphi =\bar{\Delta}_{EK\times
EL}\circ(\nabla_{EK,EL}\circ f_{EK,EL})'\circ s_{0}.$$
\end{prop}
\begin{proof}
Recall that $
\displaystyle\overline{\Delta}=\left(\sum_{j=1}^{n}\tilde{\partial}^{n+1-j}\otimes\partial_0^j
\right)\circ C(\Lambda_{EK\times EL})$ on elements of degree
$n+1$, and that $\varphi_{n} = - \varphi_{n-1}' + (\nabla
f)'s_{0}$.  The map $\varphi_{n-1}'$ is the sum of simplicial
operators whose component face and degeneracy maps all have index
at least $1$. Thus, as noted by Eilenberg and MacLane,
$\partial_{0}\varphi_{n-1}' = \varphi_{n-1}\partial_{0}$. On one
hand if $(1,x)\times(1,y)\in C_1(EK\times EL)$, then
$\varphi_0\partial_0((1,x)\times(1,y))=\varphi_0(b_0\times b_0)=0$
by construction. On the other hand if $(1,x)\times(1,y)\in
C_{k>1}(EK\times EL)$, then
$\partial_0((1,x)\times(1,y))=b_{k-1}\times b_{k-1}$ which is
degenerate, whence $\bar{\Delta}\varphi_{n-1}'$ vanishes.
\end{proof}

Thus we only need to consider the differential operator $(\nabla
f)'s_{0}$. In dimension $n>1$, $(\nabla f)'s_{0}$ is
the sum of terms of the form
$$
\begin{array}{cc}
    s_{0}^{n+1} \partial_{1} \cdots
    \partial_{n} \times s_0 &
    \hbox{ when }j=0\\
    s_{\nu_{n-j}+1} \cdots s_{\nu_1+1} s_{0} \tilde{\partial}^{n-j}
     \times s_{\mu_{j}+1} \cdots s_{\mu_1+1}
    \partial_{1}^{j-1} & \hbox{ when }0<j<n\\
    s_{0}\times s_{n} \cdots s_{1}
    \partial_{1}^{n-1} & \hbox{ when }j=n
\end{array}
$$
where $(\mu,\nu)$ is a $(j,n-j)$-shuffle, and $j$ is the running
index in the definition of the Alexander-Whitney map ($0\leq j\leq
n$). In dimension $n=1$, we have $$(\nabla
f)'s_{0}=s_0^2\partial_1\times s_0+s_0\times s_1.$$ We notice that
when $n\geq1$ the term corresponding to $j=0$ is degenerate and so
vanishes. Thus we will from now on consider $j\geq1$.

To simplify notation we
set
$s_{\mu + 1} = s_{\mu_{j}+1} \cdots s_{\mu_1+1}$,
and analogously for $s_{\nu}$.

\begin{lem}\label{lem:uno}
Let $k$ be an integer such that $1\leq k\leq n$, and let
$(\mu,\nu)$ be a $(j,n-j)$-shuffle different from the one
characterized by $\nu_{n-j}=n-j-1$, then
\begin{eqnarray}\label{degenerate}
    \partial_{0}^{k}(s_{\nu+1}s_{0}\tilde{\partial}^{n-j}
    \times s_{\mu+1}\partial_{1}^{j-1})
\end{eqnarray}
is degenerate on simplices $(1,x)\times(1,y)\in(EK \times EL)_n$.
\end{lem}
\begin{proof}
We use the simplicial identities to move one $\partial_{0}$ past
the second factor $s_{\mu+1}\partial_{1}^{j-1}$, i.e.,
\[
    \partial_{0} s_{\mu+1} \partial_{1}^{j-1} = s_{\mu}
    \partial_{0} \partial_{1}^{j-1} = s_{\mu} \partial_{0}^{j}
\]
so the second factor is the basepoint $b_{n-k+1}$ and is
degenerate since $n-k+1\neq0$ by assumption. It therefore suffices
to show that the first factor in (\ref{degenerate}) is degenerate.
Because of the simplicial identities we have
\begin{eqnarray}\label{degenerate2}
    \partial_{0}^{k}s_{\nu+1}s_{0} \tilde{\partial}^{n-j} =
    \partial_{0}^{k-1} s_{\nu} \partial_{0}s_{0}
    \tilde{\partial}^{n-j} =\partial_{0}^{k-1} s_{\nu}
    \tilde{\partial}^{n-j}.
\end{eqnarray}
Notice that if $k=1$, then (\ref{degenerate2}) is degenerate. Let
$k\geq2$, so that $n\geq2$. There exists an integer $0\leq l< n-j$
such that if $l=0$, then $\nu_1\geq1$, or if $l>0$, then
$\nu_l=l-1$ and $\nu_{l+1}>l$, so that (\ref{degenerate2}) becomes
\begin{eqnarray*}
    \partial_{0}^{k-1}s_{\nu}\tilde{\partial}^{n-j}
        & = & \partial_{0}^{k-1} s_{\nu_{n-j}} \cdots s_{\nu_{l+1}}
        s_{0}^l \tilde{\partial}^{n-j} \\
        & = &
        \left\{\begin{array}{cc}
        \partial_{0}^{k-l-1}s_{\nu_{n-j}-l} \cdots
        s_{\nu_{l+1}-l} \tilde{\partial}^{n-j} & \hbox{ if } l<k-1\\
        s_{\nu_{n-j}-k+1} \cdots
        s_{\nu_{l+1}-k+1}s_{0}^{l-k+1} \tilde{\partial}^{n-j} & \hbox{ if } l\geq k-1\\
        \end{array}\right..
\end{eqnarray*}
If $l<k-1$, then at least one $\partial_0$ can slide through
without interference to the right of the expression, i.e.,
$$
\partial_{0}^{k-l-1}s_{\nu_{n-j}-l} \cdots
        s_{\nu_{l+1}-l} \tilde{\partial}^{n-j}=
\partial_{0}^{k-l-2}s_{\nu_{n-j}-l-1} \cdots
        s_{\nu_{l+1}-l-1} \tilde{\partial}^{n-j}\partial_0.
$$
Applying this operator to $(1,x)$ we obtain the basepoint $b_{n+1-k}$, which
is degenerate since $k\leq n$. If $l\geq k-1$,
then the expression is clearly degenerate since $l<n-j$.
\end{proof}
Using Lemma~\ref{lem:uno}, we establish the following equality.
\begin{prop}\label{prop:tres}
On simplices of the form $(1,x)\times(1,y)$, we have
\[
    \bar{\Delta}\varphi\bigl((1,x) \times (1,y)\bigr)
    = \sum_{j=1}^{n}(-1)^{j(n-1)}\left( b_{n-j+1} \times
    \partial_{1}^{j-1} (1,y) \right) \otimes \left(
    \tilde{\partial}^{n-j}(1,x) \times b_{j} \right).
\]
\end{prop}
\begin{proof}
Notice that the degree of the image under $\varphi$ of
$(1,x)\times(1,y)$ is $n+1$. Thus
\begin{eqnarray}\label{tensor}
\overline{\Delta}=\left(\sum_{k=1}^n\tilde{\partial}^{n+1-k}\otimes\partial_0^k\right)
\circ C(\Lambda_{EK\times EL}),
\end{eqnarray}
where
$\tilde{\partial}^{n+1-k}=\partial_{k+1}\ldots\partial_{n+1}$. By
Lemma~\ref{lem:uno} we only need to consider the cases where the
$(j,n-j)$-shuffle $(\mu,\nu)$ in (\ref{degenerate}) is
characterized by $\nu_{n-j}=n-j-1$. Notice that the signature of
that shuffle is $\epsilon(\mu)\equiv j(n-1)\hbox{ mod }2$. Thus on
right side of the tensor in (\ref{tensor}), we have
\begin{eqnarray*}
\partial_0^k\big(s_{\nu+1}s_{0}\tilde{\partial}^{n-j}(1,x)
    &\times&
    s_{\mu+1}\partial_{1}^{j-1}(1,y)\big)\\
    &=&\partial_0^{k-1}s_{\nu}\tilde{\partial}^{n-j}(1,x)
    \times\partial_0^{k-1}s_{\mu}\partial_0^{j}(1,y)\\
    &=&\partial_0^{k-1}s_{n-j-1}\ldots s_0\tilde{\partial}^{n-j}(1,x)
    \times b_{n-k+1}\\
    &=&\partial_0^{k-1}(s_0)^{n-j}\tilde{\partial}^{n-j}(1,x)
    \times b_{n-k+1}
\end{eqnarray*}
which equals $\tilde{\partial}^{n-j}(1,x)
    \times b_{j}$, if $k-1\neq n-j$, or else, it is degenerate.
On the left side of the tensor in (\ref{tensor}), assuming that
$k-1=n-j$ and $\nu_{n-j}=n-j-1$, we have
$\mu_1=n-j,\,...,\,\mu_j=n-1$, and so
$\tilde{\partial}^{n+1-k}s_{\mu+1}=id$. Hence
\begin{eqnarray*}
\tilde{\partial}^{n+1-k}\big(s_{\nu+1}s_{0}\tilde{\partial}^{n-j}(1,x)
    &\times&
    s_{\mu+1}\partial_{1}^{j-1}(1,y)\big)\\
    &=&\partial_{n-j+2}\ldots\partial_{n+1}\big(s_{\nu+1}s_{0}\tilde{\partial}^{n-j}(1,x)
    \times
    s_{\mu+1}\partial_{1}^{j-1}(1,y)\big)\\
    &=&\partial_{n-j+2}\ldots\partial_{n+1}\big(s_{0}^{n-j+1}\tilde{\partial}^{n-j}(1,x)
    \times
    s_{\mu+1}\partial_{1}^{j-1}(1,y)\big)\\
    &=&s_{0}^{n-j+1}\partial_{1}\ldots\partial_{n}(1,x)
    \times\partial_{1}^{j-1}(1,y)\\
    &=&b_{n-j+1}\times\partial_{1}^{j-1}(1,y).
\end{eqnarray*}
The desired result is obtained by extending linearly.
\end{proof}

We now tackle $\displaystyle \prod_{k\geq1}F_k$ associated to the
E-Z data of equation (\ref{EMEZ}) on
page \pageref{EMEZ}. The first thing to notice is that the $F_k$ are
of degree $-1$, and hence vanish on elements of degree $0$.

\begin{lem}\label{lem:quatro}
If $m \geq 2$, then $F_{m}(b_{n} \times (1,y)) =
F_{m}((1,x) \times b_{n}) = 0$.
\end{lem}
\begin{proof}
The argument of Proposition~\ref{prop:form-diag} shows that
\[
    \bar{\Delta}\varphi (b_{n} \times (1,y)) = \bar{\Delta}(\nabla
    f)' s_{0}(b_{n} \times (1,y)).
\]
We can go further and show that in fact
$\bar{\Delta}\varphi (b_{n} \times (1,y)) = 0$.
Indeed, as mentioned just after Proposition~\ref{prop:form-diag}, $(\nabla f)'s_0$ is
the sum of terms
\[
    s_{\nu+1}s_{0}\tilde{\partial}^{n-j} \times
    s_{\mu+1}\partial_{1}^{j-1},
\]
where $j \geq 1$ and $(\mu,\nu)$ is a $(j,n-j)$-shuffle.  Applying
such a term to $b_{n} \times (1,y)$, we obtain
\[
    s_{\nu+1}s_{0}\tilde{\partial}^{n-j}b_{n} \times
    s_{\mu+1}\partial_{1}^{j-1}(1,y)
    = b_{n+1} \times s_{\mu+1}\partial_{1}^{j-1}(1,y).
\]
Since $j \geq 1$, the right factor is degenerate.  Since the left
factor is the basepoint, the whole expression is degenerate.
Therefore $\bar{\Delta}\varphi(b_{n} \times (1,y)) = 0$, from
which it follows that $F_{m}(b_{n} \times (1,y)) = 0$ for $m \geq
2$.  For $(1,x) \times b_{n}$ the argument is symmetric.
\end{proof}

\begin{prop}\label{prop:qF-suspension}
$\displaystyle\prod_{k\geq1}F_k=F_1\oplus F_2$ on simplices of the form
$(1,x)\times(1,y)$. Moreover,
\begin{eqnarray*}
q\circ(F_1\oplus F_2)\bigl((1,x)\!\!&\times&\!\!(1,y)\bigr)\,= \\
{[(1,x)]}\otimes[{\,}] & + & \left(\sum_{j=1}^{n}\,\,[\tilde{\partial}^{n-j}(1,x)]
\otimes [\partial_{1}^{j-1}(1,y)]\right)\quad+\quad[\,]\otimes[(1,y)].
\end{eqnarray*}
\end{prop}
\begin{proof}
A direct calculation shows that $qF_{1}((1,x) \times (1,y)) =
[(1,x)]\otimes[\,] + [\,]\otimes[(1,y)]$. Using
Proposition~\ref{prop:tres} and the fact that $q$ is a graded
algebra morphism, yields
\begin{eqnarray*}
    q\!\!&\circ&\!\!F_{2}((1,x)\times(1,y))\\
    &=&-\,q\circ F_1\otimes F_1
    \left(\sum_{j=1}^{n}(-1)^{j(n-1)}\left( b_{n-j+1} \times
    \partial_{1}^{j-1} (1,y) \right) \otimes \left(
    \tilde{\partial}^{n-j}(1,x) \times b_{j} \right)\right)\\
    &=&q\sum_{j=1}^{n}(-1)^{(n-j)+j(n-1)}
    \left([b_0\otimes\partial_1^{j-1}(1,y)\,|\,\tilde{\partial}^{n-j}(1,x)\otimes
    b_0]\right)\\
    &=&\sum_{j=1}^{n}\,(-1)^{j(n-j)+j(n-1)}\,
    [\tilde{\partial}^{n-j}(1,x)]\otimes
    [\partial_{1}^{j-1}(1,y)],
\end{eqnarray*}
where $j(n-j)+j(n-1)\equiv j(j+1)\equiv0\hbox{ mod }2$.
It remains to show that $F_{m}((1,x) \times (1,y))=0$ if $m \geq
3$.  But, by Proposition~\ref{prop:tres},
\begin{eqnarray*}
    \!\!&F_{m}&\!\!\!((1,x)\times (1,y))\,=\,-\,\sum_{i+j=m}(F_{i} \otimes
    F_{j})\bar{\Delta}\varphi ((1,x) \times (1,y)) \\
    & = &\sum_{i+j=m}\sum_{k=1}^{n}(-1)^{k(n-1)+(n-k)} F_{i}(b_{n-k+1} \times
    \partial_{1}^{k-1}(1,y)) \otimes
    F_{j}(\tilde{\partial}^{n-k}(1,x) \times b_{k}).
\end{eqnarray*}
If $m \geq 3$ then either $i \geq 2$ or $j \geq 2$, hence
Lemma~\ref{lem:quatro} implies that the above expression vanishes.
\end{proof}

\begin{cor} The map $\xi_K=q\circ\bigl(F_2\oplus F_1\bigr)\circ
C(\Lambda)$ is a twisting cochain.
\end{cor}
\begin{proof}
We need to show that $d\xi_K+\xi_Kd=\mu(\xi_K\otimes\xi_K)\Delta$.
By Lemma \ref{lem:primitive-diagonal} and the fact that $\xi_K=0$
on dimension $0$ elements, we have that
$\mu(\xi_K\otimes\xi_K)\Delta=0$. Thus we need to show that
$\xi_K$ is a chain map of degree $-1$. But, by Proposition
\ref{prop:qF-suspension}, the same proof which shows that the
Alexander-Whitney map $f_{K,L}$ is a chain map gives the desired
result.
\end{proof}

Hence the (extended) cobar diagonal for an arbitrary
pointed simplicial set $K$ is given by
\begin{eqnarray*}
\psi_{EK}\bigl([(1,x)]\bigr)&=&\\
{[(1,x)]}\otimes[\,]&+&\left(
    \sum_{j=1}^{n}\,\,[(1,\tilde{\partial}^{n-j}x)]
    \otimes [(1,\partial_{0}^{j-1}x)]\right)\quad+\quad[\,]\otimes[(1,x)],
\end{eqnarray*}
on nondegenerate simplices $(1,x)\in EK_n$.

The next corollary follows
immediately from this formula.

\begin{cor}\label{cor:aw-cobar-diagonal-suspension}
For any pointed simplicial set $K$ the map $\alpha:C(K)
\rightarrow  \Omega C(EK)$ is a coalgebra morphism with respect to
the (extended) cobar diagonal.
\end{cor}
\begin{proof}
We first check that $\alpha$ is augmentation preserving. Let
$\epsilon_K$ and $\epsilon_{EK}$ denote the augmentations on
$C(K)$ and $C(EK)$ respectively. The augmentation on
$\Omega(C(EK))$ is $\Omega(\epsilon_{EK})$ (since $\epsilon_{EK}$
is a coalgebra map). On one hand, we have
$\Omega(\epsilon_{EK})(\alpha(k_0))=\Omega(\epsilon_{EK})([\,])=1=\epsilon_K(k_0)$.
On the other hand, for $y\in K_0\backslash\{k_0\}$,
$\Omega(\epsilon_{EK})(\alpha(y-k_0))=\Omega(\epsilon_{EK})\bigl([(1,y)]\bigr)
=0=\epsilon_K(y-k_o)$.

Secondly, we check that $\alpha\otimes\alpha\circ
\Delta_K=\psi_{EK}\circ\alpha$. Let $x\in C_{n\geq1}(K)$ be a nondegenerate simplex. Then
\begin{eqnarray*}
\alpha\!\!&\otimes&\!\!\alpha\circ\Delta_K(x)\\
&=&\alpha\otimes\alpha\Bigl(\tilde{\partial}^nx\otimes x+\left(
    \sum_{j=1}^{n-1}\,\tilde{\partial}^{n-j}x
    \otimes\partial_{0}^{j}x\right)+x\otimes\partial_0^nx\Bigr)\\
&=&\alpha\otimes\alpha\Bigl(k_0\otimes
x+(\tilde{\partial}^nx-k_0)\otimes x\\
&+&\left(
    \sum_{j=1}^{n-1}\,\tilde{\partial}^{n-j}x
    \otimes\partial_{0}^{j}x\right)+x\otimes(\partial_0^nx-k_0)+x\otimes
    k_0\Bigr)\\
&=&[\,]\otimes[(1,x)]+\alpha(\tilde{\partial}^nx-k_0)\otimes
[(1,x)]\\
&+&\left(
    \sum_{j=1}^{n-1}\,[(1,\tilde{\partial}^{n-j}x)]
    \otimes[(1,\partial_{0}^{j}x)]\right)+
    [(1,x)]\otimes\alpha(\partial_0^nx-k_0)+[(1,x)]\otimes[\,].\\
\end{eqnarray*}
As in the proof of Theorem \ref{thm:geom-model-adj}, there are
four cases to consider. We will compute the case when
$\tilde{\partial}^nx\neq k_0$ and $\partial_0^nx\neq k_0$, leaving
the three other cases to the reader (one should recall that
$(1,k_0)$ is degenerate). In our case, we have
$\alpha(\tilde{\partial}^nx-k_0)=[(1,\tilde{\partial}^nx)]$ and
$\alpha(\partial_0^nx-k_0)=[(1,\partial_0^nx)]$. Now by inspection
we have $\alpha\otimes\alpha\circ
\Delta_K(x)=\psi_{EK}(\alpha(x))$.

Finally, let $y\in K_0$. By construction, we have
$\psi_{EK}([\,])=[\,]\otimes[\,]$. So clearly, we have
$\psi_{EK}(\alpha(k_0))=(\alpha\otimes\alpha)\Delta_K(k_0)$. Let
$y\in K_0\backslash\{k_0\}$. An easy computation shows that
$\psi_{EK}(\alpha(y-k_0))=\psi_{EK}\bigl([(1,y)]\bigr)=[\,]\otimes[(1,y)]+
[(1,y)]\otimes[(1,y)]+[(1,y)]\otimes[\,]$. While
\begin{eqnarray*}
\alpha\otimes\alpha\!\!&\circ&\!\!\Delta_K(y-k_0)\\
&=&\alpha\otimes\alpha(y\otimes y-k_0\otimes k_0)\\
&=&\alpha\otimes\alpha\Bigl(\bigl((y-k_0)+k_0\bigr)\otimes\bigl((y-k_0)+k_0\bigr)-k_0\otimes
k_0\Bigr)\\
&=&\alpha\otimes\alpha\Bigl((y-k_0)\otimes(y-k_0)+k_0\otimes(y-k_0)+(y-k_0)\otimes
k_0\Bigr)\\
&=&[(1,y)]\otimes[(1,y)]+[\,]\otimes[(1,y)]+[(1,y)]\otimes[\,].
\end{eqnarray*}
Hence the result.
\end{proof}

Corollary~\ref{cor:aw-cobar-diagonal-suspension}
directly implies the chain-level Bott-Samelson theorem.

\begin{thm}\label{thm:b-s-canon}  Let $K$ be any pointed simplicial set.  There is a natural isomorphism of chain Hopf algebras
$$
\widehat\alpha :(T\widetilde{C}(K), \hat d, \widehat \Delta
_{K})\rightarrow (\Om C(EK),  \psi _{EK})
$$
where $\hat d$ and $\widehat \Delta _{K}$ denote, respectively, the derivation determined by  the differential $d$ of $C(K)$ and the algebra morphism determined by $\Delta _{K}$.
\end{thm}

Therefore two of the three morphisms in (\ref{Algtriangle}) are
Hopf algebra morphisms.

\subsection{The Szczarba twisting cochain}
In this section we will assume that $K$ is a reduced pointed
simplicial set.

Since $(\Om C(EK), \psi _{EK})$ has a particularly simple form,
it is reasonable to expect that the Szczarba map also behaves better for suspensions.
We show below that our expectations are indeed fulfilled: $\theta_{EK}$ is
a strict map of chain Hopf algebras.

In \cite {S}, Szczarba defines a natural twisting cochain
$t_{L}: C(L)\longrightarrow C(GL)$, inducing a quasi-isomorphism of chain algebras
$$\theta _{L}:\Omega C(L)\longrightarrow C(GL)$$
for any $1$-reduced simplicial set $L$. Explicitly, when $L=EK$,
$t_{EK}$ is given by
\[
    t_{EK}(1,x) = \sum_{i=1}^{(n-1)!}(-1)^{\varepsilon(i,n)} D_{0,i}^{n}
    \tau((1,x))^{-1}
\]
since $\partial_{0}(1,x)=b_{n-1}$ for $x \in K_{n-1}$. The
$D_{0,i}^{n}$ are simplicial operators defined inductively as
follows.
\[
    D_{0,1}^{1} = \mathrm{id}
\]
\[
    D_{0,i+k(n-1)!}^{n+1} = \left\{
        \begin{array}{ll}
            (D_{0,i}^{n})'s_{0}\partial_{k} & k>0   \\
            (D_{0,i}^{n})' & k=0.
        \end{array}
        \right.
\]
The operator $D_{0,1}^{n}$ is therefore the identity for $n \geq
1$. The signature function $\epsilon$ is given by
\begin{eqnarray*}
\epsilon(1,1)&=&0,\\
\epsilon(i+k(n-1)!,n+1)&=&\epsilon(i,n)+k+1\quad\hbox{(mod $2$)}\\
&&\hbox{whenever }1\leq i\leq(n-1)!\hbox{, }0\leq k\leq n-1;\\
\epsilon(i,n)&=&0\quad\hbox{otherwise}.
\end{eqnarray*}

\begin{lem}\label{lem:degenerate}
$D_{0,i}^{n}$ begins with a degeneracy for all $n \geq 1$, $2 \leq
i \leq (n-1)!$.
\end{lem}
\begin{proof}
We proceed by induction on $n$. If $n=1$, then the statement is
true vacuously.  Suppose inductively that $D_{0,i}^{n}$ begins
with a degeneracy if $2 \leq i \leq (n-1)!$. For $2 \leq i \leq
n!$, consider the operator $D_{0,i}^{n+1}$.  Write $i = k(n-1)! +
\ell$, where $1 \leq \ell \leq (n-1)!$.

If $k \geq 1$, then
\[
    D_{0,i}^{n} = (D_{0,\ell}^{n})'s_{0}\partial_{k} = s_{0}
    D_{0,\ell}^{n} \partial_{k}
\]
since $D_{0,\ell}^{n}$ contains no $\partial_{0}$ (\cite{S} Lemmas
1.2 and 3.1).

If $k=0$, then $2 \leq i \leq (n-1)!$, so
\[
    D_{0,i}^{n+1} = (D_{0,i}^{n}),'
\]
which starts with a degeneracy by the inductive hypothesis.
\end{proof}

An easy induction shows that $\epsilon(1,n)=(-1)^{n+1}$. The formula for Szczarba's twisting
cochain thus becomes
\[
    t_{EK}(1,x) = (-1)^{n+1}\bigl(\tau((1,x))\bigr)^{-1}.
\]

Now a straightforward calculation proves the following theorem.

\begin{thm}\label{thm:Hopf-suspension}
The Szczarba equivalence
\[
    \theta _{EK}:(\Omega C(EK),\psi_{EK}) \xrightarrow{\simeq} C(GEK)
\]
is comultiplicative, and therefore a Hopf algebra quasi-isomorphism.
\end{thm}

\begin{rmk}
We now have two natural chain Hopf algebra quasi-isomorphisms
$\Omega C(EK)\rightarrow C(GEK)$ when $K$ is reduced. The first is
$\theta_{EK}$. The second is the composite
$$
\Omega(CEK)\xrightarrow{\gamma}C(G^+EK)\rightarrow C(GEK)
$$
where the second arrow is induced by the obvious inclusion. Note
that these two morphisms are not in general homotopic because the identity and the
inversion maps are not homotopic as simplicial maps in $GEK$.
\end{rmk}

\appendix

\section{The Milgram equivalence as a natural SDR}

Let $A$ and $B$ be coaugmented chain coalgebras. Recall that the
Milgram map $q:\Omega(A\otimes B)\rightarrow\Omega A\otimes\Omega
B$ on page \pageref{Milgram} is a quasi-isomorphism if $A$ and $B$
are 1-connected. The purpose of this Appendix is to prove the following
generalization.

\begin{thm}\label{thm:milgram-sdr}
Let $A$ and $B$ be coaugmented chain coalgebras.
There exists a strong deformation retract of chain complexes,
\[
    \Omega A \otimes \Omega B
    \sdr{\sigma}{q}
    \Omega (A \otimes B)\circlearrowleft h.
\]
In particular, $q$ is a chain homotopy equivalence.
\end{thm}

\subsection{Notation and definitions}

Throughout, we work over an arbitrary commutative ground ring $R$.

Let $(A,\Delta,\varepsilon,\eta)$ and $(B,\Delta,\varepsilon,\eta)$ be coaugmented chain
coalgebras. We adapt Sweedler's notation for the diagonal by
writing
\[
    \Delta a = a \otimes 1 + a_{1} \otimes a_{2} + 1 \otimes a.
\]
Note that we only use the notation for the reduced diagonal, and we
suppress the summation sign.
In $\Omega(A \otimes B)$, we simplify further by omitting the tensor
product symbol and the unit.  Thus $[ab]=s^{-1}(a \otimes b)$,
$[a] = s^{-1}(a \otimes 1)$, and $[b] = s^{-1}(1 \otimes b)$, so
$d[a] = -[da]+(-1)^{\deg a_{1}}[a_{1}|a_{2}]$.

Let $i_{A} = A \otimes \eta_{B} : A \cong A \otimes R
\rightarrow A \otimes B$ and $i_{B} = \eta_{A} \otimes B : B \cong
R \otimes B \rightarrow A \otimes B$.  Let $\sigma$ be the composite
\[
    \Omega A \otimes \Omega B
        \xrightarrow{\Omega i_{A} \otimes \Omega i_{B}}
            \bigl(\Omega (A \otimes B)\bigr)^{\otimes 2}
                \xrightarrow{\text{mult.}} \Omega (A \otimes B).
\]
A calculation shows that $q\sigma$ is the identity.  We show
that $q$ and $\sigma$ fit into a natural strong deformation
retract.  As a result, the Milgram map is a homotopy equivalence for
any pair of coaugmented coalgebras, not just those that are
1-connected.

\subsection{Definition of the homotopy}

In this section we construct a natural map of degree $+1$, $h :
\Omega(A \otimes B) \rightarrow \Omega(A \otimes B)$.  In the next
section we show that $h:\sigma q \simeq 1$.

Let $w \in \Omega(A \otimes B)$.  Write $w = [\omega|\beta]$,
where $\omega$ is a word that does not end in an element of $B$
and $\beta$ is possibly empty.  Define $\sharp w$ to be the total
number of letters in $\omega$ from $B$.

Define $h[\,]=0$.  If $w$ is not the empty word but ends in $ab$,
then set $h(w)=0$.

Consider $w = [ b | a_{1} | \cdots | a_{m} ]$. To define $h(w)$,
we take the iterated diagonal of $b$ and distribute the factors
among the $a_{i}$'s. We make an equivalent definition that is
easier to work with inductively. Suppose $\bar{\Delta}b = b_{1}
\otimes b_{2}$ (again, summation is understood).  Define $h[b|a] =
-(-1)^{(\deg a+1)\deg b}[ab]$.
Suppose inductively that we have defined $h[b|\alpha]$ for
any $b \in B$ and any $\alpha = [a_{1} | \cdots | a_{m-1}]$.
Define
\begin{eqnarray*}
h[ b | \alpha | a ]&=& -(-1)^{\deg b(\deg\alpha+\deg a+1)}[ \alpha | ab ]
            + ( h [ b | \alpha ] ) [ a ]\\
&&-(-1)^{\deg b_2(\deg\alpha+\deg a+1)}( h [ b_{1} | \alpha ] )
                    [ ab_{2} ].
\end{eqnarray*}

\begin{ex}
A calculation of $h[b|a|a'|a'']$, signs suppressed.
\begin{eqnarray*}
    h[b|a|a'|a'']
    & = & [a|a'|a''b] + (h[b|a|a'])[a'']
        + (h[b_{1}|a|a'])[a''b_{2}] \\
    & = & [a|a'|a''b] + [a|a'b|a'']
        + [ab|a'|a'']  \\
    &   & + [ab_{1}|a'b_{2}|a''] + [a|a'b_{1}|a''b_{2}]
        + [ab_{1}|a'|a''b_{2}] \\
    &   & + [ab_{11}|a'b_{12}|a''b_{2}].
\end{eqnarray*}
\end{ex}

Suppose that we have constructed $h$ for all words $u$ such that
$\sharp u < n$. Let $w$ be a word with $\sharp w = n$. We may
suppose that $w$ does not end in $ab$ (else we set $h(w)=0$). Thus
we may write $w = [\zeta|\alpha|\beta]$, where one of $\alpha$ or
$\beta$ may be the empty word, but not both, and $\zeta$ is either
empty or $\zeta = [\omega|x]$ with $x=b$ or $x = ab$.  If $\alpha
= [\,]$ then $x=ab$.  If $\zeta = [\,]$ or if $x = ab$ then set
$h(w)=0$. Otherwise, set
\[
    h [ \omega | b | \alpha | \beta ]
            = (-1)^{\deg\omega}[ \omega ] ( h [ b | \alpha ] ) [ \beta ]
            + (-1)^{\deg\alpha(\deg b+1)}h( [ \omega | \alpha
                | b | \beta ] ).
\]
Since $\sharp[ \omega | \alpha | b | \beta ] = \sharp w - 1$, $h[
\omega | \alpha | b | \beta ]$ has been defined.

\subsection{Induction}

In this section we show that $h:\sigma q\simeq 1$. To simplify
notation we work modulo 2.

\begin{lem}\label{lem:block}
$(dh + hd)[b | a_{1} | \cdots | a_{n} ] = [ a_{1} | \cdots | a_{n}
| b ] + [b | a_{1} | \cdots | a_{n} ]$.
\end{lem}
\begin{proof}
We proceed by induction on $n$.  When $n=1$, we have
\begin{eqnarray*}
    (dh + hd)[b|a]
    & = & d[ab] + h[db|a] + h[b|da] + h[b_{1}|b_{2}|a] +
    h[b|a_{1}|a_{2}] \\
    & = & [(da)b] + [a(db)] + [a|b] + [b|a] + [ab_{1}|b_{2}] +
        [b_{1}|ab_{2}]  \\
    & & + [a_{1}b|a_{2}] + [a_{1}|a_{2}b] +
        [a_{1}b_{1}|a_{2}b_{2}] + [a(db)] + [(da)b] \\
    & & + [b_{1}|ab_{2}] + [ab_{1}|b_{2}]
        + [a_{1}|a_{2}b] + [a_{1}b|a_{2}] + [a_{1}b_{1}|a_{2}b_{2}] \\
    & = & [a|b] + [b|a].
\end{eqnarray*}
Suppose true for $n-1$.  Let $\alpha = [a_{1}|\cdots|a_{n-1}]$.
Then
\begin{eqnarray*}
    \lefteqn{(dh+hd)[b|\alpha|a]} & & \\
    & = & d \bigl( [\alpha | ab ] + (h[b|\alpha])[a]
        + (h[b_{1}|\alpha])[ab_{2}] \bigr) \\
    &   & + h \bigl(
        [db|\alpha|a] + [b_{1}|b_{2}|\alpha|a]
            + [b](d[\alpha])[a] + [b|\alpha](d[a])
        \bigr)  \\
    & = & (d[\alpha])[ab] + [\alpha|d(ab)] + [\alpha|a|b]
        + [\alpha|b|a] \\
    &   & + [\alpha|a_{1}b|a_{2}] + [\alpha|a_{1}|a_{2}b]
        + [\alpha|ab_{1}|b_{2}] + [\alpha|b_{1}|ab_{2}] \\
    &   & + [\alpha|a_{1}b_{1}|a_{2}b_{2}]
        + (dh[b|\alpha])[a] + (h[b|\alpha])d[a]
        + (dh[b_{1}|\alpha])[ab_{2}]    \\
    &   & + (h[b_{1}|\alpha]) \bigl( [d(ab_{2})]
        + [a|b_{2}] + [b_{2}|a] + [a_{1}b_{2}|a_{2}] \\
    &   & + [a_{1}|a_{2}b_{2}] + [ab_{21}|b_{22}]
        + [b_{21}|ab_{22}] + [a_{1}b_{21}|a_{2}b_{22}] \bigr) \\
    &   & + [\alpha|a(db)] + (h[db]\alpha)[a]
        + (h[db_{1}|\alpha])[ab_{2}] +
        (h[b_{1}|\alpha])[a(db_{2})] \\
    &   & + [b_{1}](h[b_{2}|\alpha])[a] + [b_{1}|\alpha|ab_{2}]
        + [b_{1}](h[b_{21}|\alpha])[ab_{22}]
        + [\alpha|ab_{1}|b_{2}] \\
    &   & + (h[b_{1}|\alpha])[a|b_{2}]
        + (h[b_{11}|\alpha])[ab_{12}|b_{2}]
        + (d[\alpha])[ab] + (h[b]d[\alpha])[a] \\
    &   & + (h[b_{1}]d[\alpha])[ab_{2}]
        + [\alpha|(da)b] + (h[b|\alpha])[da]
        + (h[b_{1}|\alpha])[(da)b_{2}] \\
    &   & + [\alpha|a_{1}|a_{2}b] + [\alpha|a_{1}b|a_{2}]
        + (h[b|\alpha])[a_{1}|a_{2}]
        + (h[b_{1}|\alpha])[a_{1}b_{2}|a_{2}]   \\
    &   & + [\alpha|a_{1}b_{1}|a_{2}b_{2}]
        + (h[b_{1}|\alpha])[a_{1}|a_{2}b_{2}]
        + (h[b_{11}|\alpha])[a_{1}b_{12}|a_{2}b_{2}] \\
    & = & [\alpha|a|b] + [b|\alpha|a]   \\
    &   & + \bigl(
        (dh[b|\alpha])[a] + [\alpha|b|a]
        + (h[b_{1}|\alpha])[b_{2}|a] + (h[db|\alpha])[a] \\
    &   & + [b_{1}](h[b_{2}|\alpha])[a]
        + (h[b]d[\alpha])[a] + [b|\alpha|a] \bigr) \\
    &   & + \bigl(
        (dh[b_{1}|\alpha])[ab_{2}] + [\alpha|b_{1}|ab_{2}]
        + (h[b_{1}|\alpha])[b_{21}|ab_{22}] + (h[db_{1}|\alpha])[ab_{2}] \\
    &   & + [b_{1}](h[b_{21}|\alpha])[ab_{22}]
        + (h[b_{1}]d[\alpha])[ab_{2}] + [b_{1}|\alpha|ab_{2}] \bigr) \\
    & = & (\sigma q+1)[b|\alpha|a]
\end{eqnarray*}
by the inductive hypothesis and using the fact that $b_{11}
\otimes b_{12} \otimes b_{2} = b_{1} \otimes b_{21} \otimes
b_{22}$ by coassociativity.
\end{proof}

\begin{lem}\label{lem:calc}
$h([\omega](d[b|\alpha])[\beta]) = [\omega](hd[b|\alpha])[\beta] +
h([\omega](d[\alpha|b])[\beta])$.
\end{lem}
\begin{proof}
\begin{eqnarray*}
    h\bigl([\omega](d[b|alpha])[\beta]\bigr)
        & = & h \bigl( [ \omega | db | \alpha | \beta ]
            + [\omega|b_{1}|b_{2}|\alpha|\beta] \\
        &   & + [\omega | b](d[\alpha])[\beta] \bigr) \\
        & = & [\omega](h[db|\alpha])[\beta]
            + h[\omega|\alpha|db|\beta] \\
        &   & [\omega|b_{1}](h[b_{2}|\alpha])[\beta]
            + [\omega](h[b_{1}|\alpha])[b_{2}|\beta] \\
        &   & + h[\omega|\alpha|b_{1}|b_{2}|\beta]
            + [\omega](h[b]d[\alpha])[\beta] \\
        &   & + h([\omega](d[\alpha])[b|\beta]) \\
        & = & [\omega](hd[b|\alpha])[\beta]
            + h([\omega](d[\alpha|b])[\beta]).
\end{eqnarray*}
\end{proof}

\begin{lem}\label{lem:switcheroo}
$h[b|a|\alpha] = [ab|\alpha] + [a]h[b|\alpha] +
[ab_{1}]h[b_{2}|\alpha]$.
\end{lem}
\begin{proof}
Induct on length of $\alpha$.  First,
\begin{eqnarray*}
    h[b|a|a']   & = & [a|a'b] + [ab|a'] + [ab_{1}|a'b_{2}] \\
        & = & [ab|a'] + [a]h[b|a'] + [ab_{1}]h[b_{2}|a'].
\end{eqnarray*}
Consider $\alpha = [\alpha'|a']$, and suppose the lemma holds for
$[b|a|\alpha']$.  Then
\begin{eqnarray*}
    h[b|a|\alpha]
        & = & h[b|a|\alpha'|a'] \\
        & = & [a| \alpha' | a'b] + (h[b|a|\alpha'])[a']
            + (h[b_{1}|a|\alpha'])[a'] \\
        & = & [a|\alpha'|a'b] + [ab|\alpha'|a']
            + [a](h[b|\alpha'])[a'] \\
        &   & + [ab_{1}](h[b_{2}|\alpha'])[a']
            + [ab_{1}|\alpha'|a'b_{2}] \\
        &   & + [a](h[b_{1}|\alpha'])[a'b_{2}]
            + [ab_{11}](h[b_{12}|\alpha'])[a'b_{2}] \\
        & = & [ab|\alpha] + [a] \bigl( [\alpha'|a'b]
            + (h[b|\alpha'])[a'] + (h[b_{1}|\alpha'])[a'b_{2}]
            \bigr) \\
        &   & + [ab_{1}] \bigl( (h[b_{2}|\alpha'])[a']
            + [\alpha'|a'b_{2}] + (h[b_{21}|\alpha'])[a'b_{2}]
            \bigr) \\
        & = & [ab|\alpha] + [a]h[b|\alpha] +
        [ab_{1}]h[b_{2}|\alpha].
\end{eqnarray*}
In the above calculation we used the fact that $b_{1} \otimes
b_{21} \otimes b_{22} = b_{11} \otimes b_{12} \otimes b_{2}$ by
coassociativity.
\end{proof}

\begin{thm}
$dh + hd = \sigma q + 1$.
\end{thm}
\begin{proof}
Let $w \in \Omega(A \otimes B)$. Write $w = [\zeta|\alpha|\beta]$,
where $\zeta$, $\alpha$, and $\beta$ are possibly empty words,
$\alpha$ is a word from $A$, $\beta$ is a word from $B$, and
$\zeta$ is a word that does not end in a letter from $A$. If
$\zeta = [\,]$, then $(dh+hd)[\alpha|\beta] = 0 =
(\sigma q+1)[\alpha|\beta]$. If $\zeta \neq [\,]$, then write
$\zeta = [\omega|x]$, with $x = ab$ or $x = b$.
If $x=ab$ then
\begin{eqnarray*}
    (dh+hd)[\omega|ab|\alpha|\beta]
    & = & h\bigl( (d[\omega])[ab|\alpha|\beta] +
    [\omega]d[ab|\alpha|\beta] \bigr) \\
    & = & h \bigl([\omega|d(ab)|\alpha|\beta]
        + [\omega|a|b|\alpha|\beta] \\
    &   & + [\omega|b|a|\alpha|\beta]
        + [\omega|a_{1}b|a_{2}|\alpha|\beta] \\
    &   & + [\omega|a_{1}|a_{2}b|\alpha|\beta]
        + [\omega|ab_{1}|b_{2}|\alpha|\beta] \\
    &   & + [\omega|b_{1}|ab_{2}|\alpha|\beta]
        + [\omega|a_{1}b_{1}|a_{2}b_{2}|\alpha|\beta] \bigr) \\
    & = & [\omega|a](h[b|\alpha])[\beta]
        + h[\omega|a|\alpha|b|\beta] \\
    &   & + [\omega|ab|\alpha|\beta]
        + [\omega|a](h[b|\alpha])[\beta] \\
    &   & + [\omega|ab_{1}](h[b_{2}|\alpha])[\beta]
        + h([\omega|a|\alpha|b|\beta]) \\
    &   & + [\omega|ab_{1}](h[b_{2}|\alpha])[\beta]
        + h[\omega|ab_{1}|\alpha|b_{2}|\beta] \\
    & = & [\omega|ab|\alpha|\beta]
\end{eqnarray*}
where we used Lemma~\ref{lem:switcheroo}.

Suppose that $x=b$.  We induct on $\sharp w$. The case $\sharp w =
1$, $\omega$ and $\beta$ empty was treated in
Lemma~\ref{lem:block} and the general case of $\sharp w = 1$ is an
easy calculation using Lemma~\ref{lem:calc}.

Suppose that $(dh+hd)(u) = \sigma q(u) + u$ for all words $u$
such that $\sharp u < n$.  Let $w$ be a word with $\sharp w = n$.
Write $w = [\omega | b | \alpha | \beta]$. Note that the inductive
hypothesis applies to $[b|\alpha]$ and to
$[\omega|\alpha|b|\beta]$.  Thus
\begin{eqnarray*}
    (dh+hd)(w)
    & = & d \bigl( [w](h[b|\alpha])[\beta] \bigr)
        + dh[\omega|\alpha|b|\beta] \\
    &   & + h\bigl( (d[\omega]) [b | \alpha|\beta]
        +  [\omega](d[b|\alpha])[\beta]
        + [\omega|b|\alpha]d[\beta] \bigr) \\
    & = & (d[\omega])(h[b|\alpha])[\beta]
        + [\omega](dh[b|\alpha])[\beta]
        + [\omega](h[b|\alpha])d[\beta] \\
    &   & + dh[\omega|\alpha|b|\beta]
        + (d[\omega])(h[b|\alpha])[\beta]
        + h\bigl((d[\omega])[\alpha|b|\beta] \bigr) \\
    &   & + [\omega](hd[b|\alpha])[\beta]
        + h\bigl([\omega](d[\alpha|b])[\beta] \bigr)
        + [\omega](h[b|\alpha])d[\beta] \\
    &   & + h\bigl( [\omega|\alpha|b]d[\beta] \bigr) \\
    & = & [\omega]\bigl( [\alpha|b] + [b|\alpha] \bigr) [\beta]
        + \sigma q[\omega|\alpha|b|\beta]
        + [\omega|\alpha|b|\beta] \\
    & = & (\sigma q + 1)(w)
\end{eqnarray*}
where we have used Lemma~\ref{lem:calc} and the fact that
$\sigma q[\omega|\alpha|b|\beta] = \sigma q(w)$.
\end{proof}

\subsection{Elementary properties}

\begin{prop}\label{prop:elem}
The natural homotopy $h$ satisfies
\begin{enumerate}
    \item\label{elem-1} $qh = 0$,
    \item\label{elem-2} $h\sigma = 0$, and
    \item\label{elem-3} $h^{2} = 0$.
\end{enumerate}
\end{prop}
\begin{proof}
Properties (\ref{elem-1}) and (\ref{elem-2}) follow from the
definitions.  Property (\ref{elem-3}) is proved by an induction on
$\sharp w$ and the length of the last block of letters from $A$.
\end{proof}

\end{document}